\documentclass[11pt]{amsart}
\usepackage{amsmath,amssymb,amsthm,array}
\usepackage[colorlinks=true]{hyperref}
\theoremstyle{plain}
\newtheorem{thm}{Theorem}[section]
\newtheorem{prop}[thm]{Proposition}
\newtheorem{lem}[thm]{Lemma}
\newtheorem{cor}[thm]{Corollary}
\theoremstyle{definition}
\newtheorem{defn}[thm]{Definition}
\newtheorem{exam}[thm]{Example}

\theoremstyle{remark}
\newtheorem{rem}[thm]{Remark}
\newcommand{\F}{\mathbb F}
\newcommand{\Hull}{\operatorname{Hull}}
\newcommand{\ev}{\operatorname{ev}}
\newcommand{\supp}{\operatorname{supp}}
\newcommand{\res}{\operatorname{res}}
\begin{document}
\title[Geometric Obstructions to LCD Scalings]{Geometric Obstructions and Sharp Genus Dependence in LCD Scalings}
\author{Rahim Rahmati-Asghar}
\address{Department of Mathematics, Faculty of Basic Sciences, University of Maragheh, Maragheh, Iran}
\email{rahmatiasghar.r@maragheh.ac.ir, rahmatiasghar.r@gmail.com}
\subjclass[2020]{Primary 94B27; Secondary 14H05, 11T71}
\keywords{Algebraic-geometric code, Euclidean hull, LCD code, Riemann--Roch space, geometric threshold}

\begin{abstract}
We study the least pole order of a rational function whose nonzero evaluations scale a prescribed code to a Euclidean LCD code. A residue-one differential yields an obstruction for odd-dimensional one-point AG codes in characteristic two: for $C=C_L(D,aQ)$ with $0\le a<n$ and odd dimension, the threshold is at least $n+2g-a-1$. The obstruction uses alternating Gram matrices and does not require the scaled code to be self-orthogonal. For the dimension-three codes $C_L(D,(r+1)Q)$ on Hermitian curves over $\F_{r^2}$, $r=2^e\ge4$, we obtain
\[
n+2g-r-2\le\tau_D(C;Q)\le n+2g-3.
\]
Consequently the leading genus contribution is asymptotically sharp even in dimension three. On the elliptic curve $y^2+y=x^3$ over $\F_4$, we prove the exact dimension-three threshold $6$ by an explicit multiplier. For repetition codes on the Hermitian family, the threshold equals the general interpolation bound $n+2g-1$. The upper bounds use information sets and hold for $q>3k$, and for $q>2k$ in characteristic two.
\end{abstract}
\maketitle

\section{Introduction}
Every linear code over $\F_q$ with $q>3$ is monomially equivalent to a Euclidean linear complementary dual (LCD) code \cite{CarletEtAl2018}. When the coordinates are indexed by rational points of a curve, one may ask for the least pole order of a single function realizing these multipliers.

Let $X/\F_q$ be a smooth, projective, geometrically integral curve of genus $g$, let $D=P_1+\cdots+P_n$ consist of distinct rational points, and let $Q$ be rational and outside $\supp(D)$. For a code $C\subseteq\F_q^n$, let $\tau_D(C;Q)$ be the least $m$ for which a function $f\in L(mQ)$ is nonzero at every $P_i$ and scales $C$ to an LCD code. For $q>3$, interpolation gives $\tau_D(C;Q)\le n+2g-1$. We ask when the genus contribution is forced by an obstruction to all smaller pole orders.

The main obstruction concerns odd-dimensional codes in characteristic two. Suppose a differential has divisor $(n+2g-2)Q-D$ and residue one at every $P_i$. If $C=C_L(D,aQ)$ has odd dimension, then
\[
\tau_D(C;Q)\ge n+2g-a-1\qquad(0\le a<n,\ a\le n+2g-2).
\]
Below this bound the Gram matrix of every scaled code is alternating and hence singular. This argument is stronger than requiring the entire Gram matrix to vanish.

For the Hermitian curve $y^r+y=x^{r+1}$ over $\F_{r^2}$, with $r$ a power of two, take all $n=r^3$ affine rational points and the point $Q$ at infinity. The code $C_L(D,(r+1)Q)$ has dimension three. For $r\ge4$, our lower bound and an information-set upper bound differ by only $r-1$:
\[
 n+2g-r-2\le\tau_D(C_L(D,(r+1)Q);Q)\le n+2g-3.
\]
Since $g=r(r-1)/2$, these inequalities prove $\tau_D(C;Q)=n+2g-O(\sqrt g)$ along this family. Thus even in fixed dimension three the coefficient of $g$ cannot be uniformly reduced below two. 

\subsection*{Relation to existing results}
The following table compares the frameworks of the cited papers with the quantity studied here.

\smallskip
\noindent\begin{tabular}{@{}p{0.23\textwidth}p{0.34\textwidth}p{0.36\textwidth}@{}}
\hline
Source & Established framework & Quantity studied in this paper\\\hline
Carlet et al.\ \cite{CarletEtAl2018} & Gram determinants and coordinate scalings yielding LCD equivalents & Pole cost of realizing multipliers in $L(mQ)$\\
Jin--Xing \cite{JinXing2018} & Counting admissible scalings of AG codes in even characteristic & Constraints on every low-pole multiplier\\
Mesnager--Tang--Qi \cite{MesnagerTangQi2018} & Residue multipliers and non-special-divisor LCD criteria & Minimum pole order for fixed $C,D,Q$\\
Camps Moreno et al.\ \cite{CampsMorenoLopezMatthews2024} & Divisor intersections describing AG hulls & Bounds uniform over an entire multiplier space\\
Bras-Amor\'os et al.\ \cite{BrasAmorosDuursmaHong2019} & Self-dual flags and their terminal pole orders & LCD thresholds constrained by the same filtration\\\hline
\end{tabular}
\smallskip

The coordinate-support idea in the upper bound is related to systematic-generator scaling methods in \cite{CarletEtAl2018}; the extra step is to realize a suitable perturbation through a low-pole function. The residue differential on the Hermitian curve is already used in \cite[Section V-B]{MesnagerTangQi2018}. Standard AG duality \cite{Stichtenoth2009} identifies $C_L(D,aQ)^\perp=C_L(D,(n+2g-2-a)Q)$ in this setting. In particular, the repetition threshold below is a direct consequence of this duality and Frobenius, rather than a new Hermitian duality result. The higher-dimensional argument combines the same containment with the elementary odd-order alternating-matrix obstruction.

There is also a close connection with isometry-dual flags \cite{BrasAmorosDuursmaHong2019}, where the terminal evaluation step at $n+2g-1$ and examples on Hermitian curves are already studied. This explains why the repetition threshold is closely tied to a classical terminal pole-order phenomenon. Our dimension-three obstruction concerns singularity of every scaled Gram matrix along the filtration, rather than the characterization of an isometry-dual flag.

Our contributions are the pole-order obstructions, the explicit threshold calculation, and the geometric upper bounds.

\section{The threshold and the interpolation bound}
Throughout, $q$ is a prime power; hulls and orthogonality are Euclidean. Field-size hypotheses are stated in the individual results.

For a divisor $A$, set
\[
L(A)=\{0\}\cup\{f\in\F_q(X)^\times:(f)+A\ge0\},
\qquad \ell(A)=\dim_{\F_q}L(A).
\]
If $\supp(A)\cap\supp(D)=\varnothing$, the AG evaluation code is $C_L(D,A)=\ev_D(L(A))$. The evaluation kernel is $L(A-D)$, and Riemann--Roch gives $\ell(A)=\deg A+1-g$ for $\deg A>2g-2$ \cite{Stichtenoth2009}.

\begin{defn}[Euclidean dual, hull, and LCD code]\label{def:lcd}
Let $C\subseteq\F_q^n$ be an $[n,k]$ linear code. Its Euclidean dual is
\[
C^\perp=\{v\in\F_q^n:\langle v,c\rangle=0\text{ for every }c\in C\},
\qquad \langle v,c\rangle=\sum_{i=1}^n v_ic_i.
\]
The Euclidean hull of $C$ and its hull dimension are
\[
\Hull(C)=C\cap C^\perp,\qquad h(C)=\dim_{\F_q}\Hull(C).
\]
The code $C$ is called \emph{linear complementary dual} (LCD) if $C\cap C^\perp=\{0\}$ \cite{Massey1992}. Equivalently,
\begin{equation}\label{eq:lcd-equivalence}
C\text{ is LCD}\quad\Longleftrightarrow\quad h(C)=0
\quad\Longleftrightarrow\quad\F_q^n=C\oplus C^\perp.
\end{equation}
The last equivalence follows from $\dim C^\perp=n-k$. In general,
$0\le h(C)\le\min\{k,n-k\}$; moreover, $h(C)=k$ if and only if $C$ is self-orthogonal, meaning $C\subseteq C^\perp$.
\end{defn}

Let $B$ be a full-row-rank generator matrix of $C$. The Gram-matrix characterization \cite{Massey1992} gives
\begin{equation}\label{eq:unscaled-gram}
h(C)=k-\operatorname{rank}(BB^T).
\end{equation}
Thus $C$ is LCD precisely when $BB^T$ is nonsingular. More generally, an invertible diagonal scaling $\Lambda$ has generator matrix $B\Lambda$, so
\begin{equation}\label{eq:gram}
\dim\Hull(C\Lambda)=k-\operatorname{rank}(B\Lambda^2B^T).
\end{equation}
Indeed, the injective map $u\mapsto uB\Lambda$ identifies the hull with the left kernel of $B\Lambda^2B^T$: the vector $uB\Lambda$ is orthogonal to every row of $B\Lambda$ exactly when $uB\Lambda^2B^T=0$. In particular, finding an LCD scaling is equivalent to making this Gram matrix have full rank $k$, or, equivalently, reducing the hull dimension to zero.

\begin{defn}
For any finite field $\F_q$, define
\[
\begin{split}
\tau_D(C;Q)=\min\bigl\{m\in\mathbb Z_{\ge0}:{}&\exists f\in L(mQ),\\
 &f(P_i)\ne0\ (1\le i\le n),\ \Hull(C\Lambda_f)=\{0\}\bigr\}.
\end{split}
\]
where $\Lambda_f=\operatorname{diag}(f(P_1),\ldots,f(P_n))$, with the convention $\min\varnothing=+\infty$.
\end{defn}
The curve, evaluation points, and pole point are part of the data. Although $C$ may be any linear code, the allowed multipliers are constrained by these data. Since $L(0)=\F_q$, one has $\tau_D(C;Q)=0$ if and only if $C$ is LCD.

The definition itself needs no restriction $q>3$. This hypothesis guarantees finiteness for every code through LCD equivalence and interpolation, as proved below. It is not needed for the residue obstructions. Over $\F_2$ or $\F_3$, every nonzero scalar has square one, so $B\Lambda^2B^T=BB^T$ for every invertible diagonal $\Lambda$. Consequently a non-LCD code over either field has threshold $+\infty$. Examples are the repetition codes of length two over $\F_2$ and length three over $\F_3$, evaluated at the affine points of the respective projective lines.

\begin{prop}[Interpolation bound]\label{prop:universal}
For $q>3$, the threshold is finite and $\tau_D(C;Q)\le n+2g-1$.
\end{prop}
\begin{proof}
Set $m=n+2g-1$. Both $mQ$ and $mQ-D$ have degree greater than $2g-2$, so Riemann--Roch gives $\ell(mQ)-\ell(mQ-D)=n$. Thus evaluation onto $\F_q^n$ is surjective. The LCD-equivalence theorem \cite{CarletEtAl2018}, together with invariance under coordinate permutations, supplies a diagonal LCD scaling with nonzero entries $\lambda_i$. Interpolate $f(P_i)=\lambda_i$ in $L(mQ)$. The evaluation-surjectivity argument is valid over every finite field; $q>3$ is used only to obtain the LCD scaling.
\end{proof}

\section{A higher-dimensional residue obstruction}\label{sec:sharp}
Put $s=n+2g-2$ and suppose there is a nonzero rational differential satisfying
\begin{equation}\label{eq:omega}
(\omega)=sQ-D,\qquad \res_{P_i}(\omega)=1\quad(1\le i\le n).
\end{equation}
The global residue theorem gives the standard identity
\begin{equation}\label{eq:residuesum}
\sum_i h(P_i)=0\qquad(h\in L(sQ)).
\end{equation}
Equivalently, AG duality gives $C_L(D,aQ)^\perp=C_L(D,(s-a)Q)$ \cite{Stichtenoth2009}. For comparison, we first record the usual product-space obstruction.

\begin{prop}[Standard product-space obstruction]\label{prop:ordinary}
Assume \eqref{eq:omega}, let $0\le a<n$, and put $C=C_L(D,aQ)$. For an admissible $f\in L(mQ)$, set $b=\min\{a,s-a-2m\}$. Then
\[
\dim\Hull(C\Lambda_f)\ge\ell(bQ).
\]
In particular, if $s-a\ge0$, then $\tau_D(C;Q)\ge\lfloor(s-a)/2\rfloor+1$.
\end{prop}
\begin{proof}
For $u\in L(bQ)$ and $v\in L(aQ)$, the product $f^2uv$ belongs to $L(sQ)$, so \eqref{eq:residuesum} gives $\sum_i f(P_i)^2u(P_i)v(P_i)=0$. Thus $\ev_D(u)\Lambda_f$ belongs to the hull. Evaluation on $L(bQ)$ is injective because $b\le a<n$, and diagonal scaling is invertible. If $b\ge0$, the constant functions already give a nonzero hull vector. This yields the final assertion.
\end{proof}

\begin{thm}[Odd-dimensional pole-order obstruction]\label{thm:odd}
Assume \eqref{eq:omega} over a field of characteristic two. Let $0\le a<n$, $a\le s$, and suppose $C=C_L(D,aQ)$ has odd dimension $k$. Then
\begin{equation}\label{eq:odd}
\tau_D(C;Q)\ge s-a+1=n+2g-a-1.
\end{equation}
More precisely, every admissible $f\in L(mQ)$ with $m+a\le s$ gives an odd-dimensional hull.
\end{thm}
\begin{proof}
Since $a<n$, evaluation is injective on $L(aQ)$. Choose a basis $u_1,\ldots,u_k$ of this space and form the Gram matrix $M_f$ of the scaled evaluation vectors. Its diagonal entries are
\[
(M_f)_{jj}=\sum_i (fu_j)(P_i)^2
 =\left(\sum_i(fu_j)(P_i)\right)^2=0,
\]
because $fu_j\in L((m+a)Q)\subseteq L(sQ)$. The matrix $M_f$ is symmetric with zero diagonal, hence alternating. An alternating matrix has even rank: a nonzero off-diagonal entry splits off a nondegenerate alternating block of order two, and induction applies to its orthogonal complement. Since $k$ is odd, $k-\operatorname{rank}M_f$ is a positive odd integer. Equation~\eqref{eq:gram} proves both assertions.
\end{proof}

Theorem \ref{thm:odd} holds over $\F_2$ as well. The hypothesis $q>3$ enters Corollary \ref{thm:residue} through the interpolation upper bound.

\begin{rem}
The improvement comes from parity; the Gram matrix need not vanish. An alternating Gram matrix can have nonzero off-diagonal entries. Odd dimension is essential: an even-order alternating matrix can be nonsingular.
\end{rem}

\begin{cor}[Exact repetition threshold]\label{thm:residue}
Under \eqref{eq:omega}, if $q>3$ has characteristic two and $s\ge0$, then
\[
\tau_D(\langle\mathbf1\rangle;Q)=n+2g-1.
\]
\end{cor}
\begin{proof}
Apply Theorem~\ref{thm:odd} with $a=0$, $k=1$, and combine it with Proposition~\ref{prop:universal}.
\end{proof}

\begin{thm}[Exact Hermitian threshold]\label{thm:hermitian}
Let $r=2^e$ with $e\ge1$, and let $X/\F_{r^2}$ be the smooth projective Hermitian curve with affine equation
\[
y^r+y=x^{r+1}.
\]
Let $Q$ be its unique point at infinity, and let $D$ be the sum of all its affine $\F_{r^2}$-rational points. Then $n=r^3$, $g=r(r-1)/2$, and
\begin{equation}\label{eq:sharp}
\tau_D(C_L(D,0);Q)=r^3+r(r-1)-1=n+2g-1.
\end{equation}
\end{thm}
\begin{proof}
The standard Hermitian-curve formulas give the genus, the unique point at infinity, and $(x)_\infty=rQ$ \cite{Stichtenoth2009}. For each $a\in\F_{r^2}$, $a^{r+1}\in\F_r$, and the trace equation $y^r+y=a^{r+1}$ has exactly $r$ solutions in $\F_{r^2}$. Hence there are $r^3$ affine rational points.

Put $h=x^{r^2}-x$. At every affine point, the partial derivative of the defining equation with respect to $y$ is one. Thus $x-a$ is a local parameter at each point above $x=a$, and $dx$ is regular and nonvanishing on the affine curve. Since a differential has degree $2g-2$,
\[
(dx)=(2g-2)Q.
\]
The zeros of $h$ are precisely the affine rational points, all simple; its unique pole has order $r^3$. Therefore
\[
(h)=D-nQ,\qquad
\omega=\frac{dh}{h}=\frac{-dx}{x^{r^2}-x}
\]
satisfies $(\omega)=(n+2g-2)Q-D$. At each simple zero of $h$, its logarithmic differential $dh/h$ has residue one. Corollary~\ref{thm:residue} applies because $r^2\ge4$.
\end{proof}

\begin{rem}
The preceding corollary and theorem establish equality in the interpolation bound. The differential itself is classical \cite[Section V-B]{MesnagerTangQi2018}. Here ``Hermitian'' describes the curve; every hull in this paper is Euclidean.
\end{rem}

\section{An information-set upper bound}\label{sec:main}
An information set is a set $I$ of $k$ coordinates for which the square submatrix $B_I$ is invertible. If $v$ has support exactly $I$, then
\begin{equation}\label{eq:supported}
\det(B\operatorname{diag}(v_i^2)B^T)
=\det(B_I)^2\prod_{i\in I}v_i^2\ne0.
\end{equation}
Indeed, only the columns in $I$ contribute. Equivalently, this is the single surviving term in the Cauchy--Binet expansion \cite{HornJohnson2013}.

\begin{lem}\label{lem:affine}
Suppose $1\le k<n$, $C$ is not LCD, and $v\in\F_q^n$ has support exactly an information set $I$. If $q>3k$, or if the characteristic is two and $q>2k$, some $t\in\F_q$ makes $v+t\mathbf1$ nonzero in every coordinate and gives an LCD scaling of $C$.
\end{lem}
\begin{proof}
Set
\[
p(T)=\det(B\operatorname{diag}((v_i+T)^2)B^T).
\]
By \eqref{eq:supported}, $p(0)\ne0$. Its degree is at most $2k$, and the coefficient of $T^{2k}$ is $\det(BB^T)=0$ since $C$ is not LCD. Hence it has at most $2k-1$ roots. Coordinate nonvanishing excludes $t=0$ and at most $k$ further values $-v_i$, $i\in I$. At most $3k$ values are excluded altogether.

In characteristic two, Cauchy--Binet gives $p(T)=h(T)^2$, where
\[
h(T)=\sum_{|J|=k}\det(B_J)\prod_{j\in J}(v_j+T).
\]
The leading coefficient of $h$ squares to $\det(BB^T)=0$, so $\deg h\le k-1$. There are then at most $k-1$ determinant roots and at most $k+1$ inadmissible values. Thus $q>2k$ suffices. Equation~\eqref{eq:gram} completes the proof.
\end{proof}

\begin{thm}[Information-set bound]\label{thm:main}
For any $[n,k]$ linear code with $1\le k<n$, if $q>3k$, or if the characteristic is two and $q>2k$, then
\begin{equation}\label{eq:main}
\tau_D(C;Q)\le n-k+2g.
\end{equation}
In particular, this holds for every AG evaluation code of dimension $k$.
\end{thm}
\begin{proof}
If $C$ is LCD, use $f=1$. Otherwise fix an information set $I$, put $D_{I^c}=\sum_{j\notin I}P_j$, and set $m=n-k+2g$ and $A=mQ-D_{I^c}$. Then $\deg A=2g$ and $\deg(A-P_i)=2g-1$ for $i\in I$. Riemann--Roch gives
\[
\ell(A)-\ell(A-P_i)=1.
\]
Thus each evaluation functional $L(A)\to\F_q$, $f\mapsto f(P_i)$, is nonzero. Their $k$ kernels cover at most $kq^{\ell(A)-1}<q^{\ell(A)}$ elements. Choose $f$ outside their union. Its evaluation vector has support exactly $I$. Lemma~\ref{lem:affine} supplies $t$ such that $f+t\in L(mQ)$ is admissible and gives an LCD scaling.
\end{proof}

\begin{rem}[Sufficient field-size conditions]\label{rem:field-size}
The conditions $q>3k$ and, in characteristic two, $q>2k$ are sufficient conditions from the root counts in Lemma~\ref{lem:affine}; no optimality or necessity is asserted. The count bounds the determinant roots and coordinate exclusions separately, allowing both overlap and repeated excluded values, so particular instances can work over smaller fields. In particular, an already LCD code has threshold zero regardless of these conditions. The exact example in Theorem~\ref{thm:elliptic} succeeds with $q=4$ and $k=3$, outside both sufficient ranges. Example~\ref{exam:parity} shows attainment of the pole-order bound, not sharpness of its field-size hypothesis. Uniformly weakening these hypotheses for arbitrary non-LCD codes remains open here.
\end{rem}

\begin{cor}[Divisor-sensitive version]\label{cor:mu}
For an information set $I$, let $\mu_I$ be the least nonnegative integer $m$ for which
\[
\ell(mQ-D_{I^c})>\ell(mQ-D_{I^c}-P_i)\qquad(i\in I).
\]
Under either field-size hypothesis of Theorem~\ref{thm:main}, with $k<n$,
\[
\tau_D(C;Q)\le\min_{I\text{ information set}}\mu_I\le n-k+2g.
\]
\end{cor}
\begin{proof}
The dimension inequalities make all the relevant evaluation functionals nonzero. If $C$ is not LCD, the same union-of-kernels and constant-perturbation argument applies at $m=\mu_I$. If $C$ is LCD, its threshold is zero. The final inequality follows from the degree calculation in Theorem~\ref{thm:main}.
\end{proof}

\begin{exam}[Attainment in genus zero]\label{exam:parity}
Take five distinct affine rational points on $\mathbb P^1_{\F_{25}}$, let $D$ be their sum, and set $Q=P_\infty$. The single-parity-check code
\[
C=\{(c_1,\ldots,c_5)\in\F_{25}^5:\textstyle\sum_i c_i=0\}
\]
has dimension four and $C^\perp=\langle\mathbf1\rangle\subseteq C$ because the characteristic is five. It is not LCD, so its threshold is positive. Since $25>3\cdot4$, Theorem~\ref{thm:main} gives
\[
\tau_D(C;Q)=1=n-k+2g.
\]
This example proves attainment of the information-set bound in genus zero, separately from the genus dependence established in Section~\ref{sec:sharp}.
\end{exam}

\section{Dimension three: sharp genus dependence and an exact threshold}
\begin{thm}[A dimension-three Hermitian family]\label{thm:three}
On the Hermitian curve of Theorem~\ref{thm:hermitian}, let $a=r+1$ and $C=C_L(D,aQ)$. Then $\dim C=3$. For $r\ge4$,
\begin{equation}\label{eq:three}
n+2g-r-2\le\tau_D(C;Q)\le n+2g-3.
\end{equation}
In particular,
\[
\lim_{r=2^e\to\infty}\frac{\tau_D(C;Q)-n}{g}=2.
\]
Consequently, there are no constants $\varepsilon>0$ and
$K\in\mathbb R$, independent of $r$, such that
\[
\tau_D(C;Q)\le n+(2-\varepsilon)g+K
\]
for all sufficiently large powers of two $r$.
\end{thm}
\begin{proof}
The standard one-point Hermitian basis consists of $x^iy^j$ with $i\ge0$, $0\le j<r$, and pole order $ri+(r+1)j$ \cite{Stichtenoth2009}. At $a=r+1$, it gives precisely $1,x,y$. Evaluation is injective since $a<n$. The differential constructed in Theorem~\ref{thm:hermitian} satisfies \eqref{eq:omega}, so Theorem~\ref{thm:odd} gives the lower bound. Since $q=r^2>2\cdot3$ for $r\ge4$, Theorem~\ref{thm:main} gives the upper bound. 

For the final assertion, suppose that constants $\varepsilon>0$
and $K$, independent of $r$, satisfy
\[
\tau_D(C;Q)\le n+(2-\varepsilon)g+K
\]
throughout the family. Dividing by $g$ gives
\[
\frac{\tau_D(C;Q)-n}{g}
\le 2-\varepsilon+\frac{K}{g}.
\]
Since $g=r(r-1)/2\to\infty$, taking limits yields
$2\le 2-\varepsilon$, contradicting $\varepsilon>0$.
\end{proof}

\begin{thm}[An exact positive-genus threshold in dimension three]\label{thm:elliptic}
Let $X/\F_4$ be $y^2+y=x^3$, let $D$ contain all eight affine rational points, and let $Q$ be the point at infinity. Then
\[
\tau_D(C_L(D,3Q);Q)=6.
\]
An admissible multiplier attaining this value is $f=1+x^3+xy$.
\end{thm}
\begin{proof}
Here $g=1$, $s=8$, and $C_L(D,3Q)$ has basis $1,x,y$. Theorem~\ref{thm:odd} gives $\tau_D(C;Q)\ge8-3+1=6$. Since $x$ and $y$ have pole orders two and three, respectively, $f\in L(6Q)$.

Write $\F_4=\{0,1,\alpha,\alpha+1\}$ with $\alpha^2+\alpha+1=0$. Above $x=0$ the two points have $y=0,1$, and $f=1$. Above each nonzero $x$, one has $x^3=1$ and $y\in\{\alpha,\alpha+1\}$, so $f=xy\ne0$. Thus $f$ is admissible. The complete evaluation data are
\[
\begin{array}{c|cccccccc}
x&0&0&1&1&\alpha&\alpha&\alpha+1&\alpha+1\\
y&0&1&\alpha&\alpha+1&\alpha&\alpha+1&\alpha&\alpha+1\\
f&1&1&\alpha&\alpha+1&\alpha+1&1&1&\alpha
\end{array}
\]
For the generator matrix $B$ formed by evaluating $1,x,y$, these values give
\[
B\Lambda_f^2B^T=
\begin{pmatrix}0&1&1\\1&0&0\\1&0&1\end{pmatrix},
\qquad \det(B\Lambda_f^2B^T)=1.
\]
This can also be checked fiberwise using $\alpha+(\alpha+1)=1$ and $\sum_{x\in\F_4^*}x=0$. Equation~\eqref{eq:gram} proves the upper bound six. The information-set field-size hypothesis is not used for this exact example.
\end{proof}

\begin{exam}
For $r=4$, Theorem~\ref{thm:three} gives $70\le\tau_D(C_L(D,5Q);Q)\le73$ over $\F_{16}$, with $n=64$ and $g=6$. These are bounds, not an asserted exact threshold. For comparison, the repetition code on the same curve has exact threshold $75$ by Theorem~\ref{thm:hermitian}.
\end{exam}

\section{The projective line and the role of characteristic}
Let $X=\mathbb P^1_{\F_q}$, $Q=P_\infty$, $D=\sum_{a\in\F_q}P_a$, and $R=\langle\mathbf1\rangle$. Then $L(mQ)$ consists of polynomials of degree at most $m$.

\begin{prop}\label{prop:repetition}
For $q>3$,
\[
\tau_D(R;P_\infty)=
\begin{cases}
(q-1)/2,&q\text{ odd},\\
q-1,&q\text{ even}.
\end{cases}
\]
\end{prop}
\begin{proof}
Suppose first that $q$ is odd and put $d=(q-1)/2$. The finite-field power sums satisfy $\sum_a a^j=0$ for $0\le j<q-1$ and $\sum_a a^{q-1}=-1$ \cite{LidlNiederreiter1997}. If $\deg f<d$, then $\sum_a f(a)^2=0$. Conversely choose $c\notin\{0,1,-1\}$ and $f=x^d+c$. Its evaluations are nonzero because $a^d\in\{1,-1\}$ for $a\ne0$, and
\[
\sum_a f(a)^2=\sum_a(a^{q-1}+2ca^d+c^2)=-1\ne0.
\]
This proves the odd case. If $q$ is even, every polynomial of degree below $q-1$ has zero evaluation sum, hence zero sum of squared evaluations. This gives the lower bound $q-1$, and Proposition~\ref{prop:universal} gives the upper bound. Equivalently, the even case follows from Corollary~\ref{thm:residue} with $\omega=d(x^q-x)/(x^q-x)$ and $s=q-2$.
\end{proof}

\section{Application to constrained LCD conversion}\label{sec:application}
Orthogonal direct-sum masking uses an LCD code to separate data and random-mask components \cite{CarletGuilley2016}; the associated orthogonal recovery requires a nonsingular Gram matrix. The results here provide \emph{feasibility certificates} when conversion of $C_L(D,aQ)$ is restricted to multipliers $\ev_D(f)$ with $f\in L(mQ)$. Under the hypotheses of Theorem~\ref{thm:odd}, every budget $m\le n+2g-a-2$ is impossible, whereas Theorem~\ref{thm:main} gives a sufficient budget under its field-size conditions. Entire multiplier spaces can therefore be excluded before candidate enumeration.

For $C=C_L(D,5Q)$ on $y^4+y=x^5$ over $\F_{16}$, Theorem~\ref{thm:three} yields
\[
\begin{array}{c|l}
\text{Permitted budget}&\text{Conclusion}\\\hline
m\le69&\text{No admissible multiplier can give an LCD code.}\\
70\le m\le72&\text{Existence is not decided by the present bounds.}\\
m\ge73&\text{An admissible LCD multiplier exists.}
\end{array}
\]
The last row uses $L(73Q)\subseteq L(mQ)$. These conclusions apply to the restricted geometric conversion problem. Pole order is an algebraic budget, not a hardware cost.

\section{Scope and further questions}
The odd-dimensional obstruction shows that low-pole scalings can remain non-LCD after the usual self-orthogonality obstruction has ceased to apply. The dimension-three Hermitian family forces a pole order $n+2g-O(\sqrt g)$, while the elliptic example gives a fully explicit exact threshold beyond dimension one. Neither the code parameters nor the existence of LCD equivalents is improved: diagonal scaling preserves length, dimension, and minimum distance.

The bounds in Theorem~\ref{thm:three} differ by $r-1$; their equality is not asserted for $r\ge4$. Determining these thresholds exactly is a concrete next question. Other unresolved issues are the even-dimensional case, where alternation need not cause singularity, and odd characteristic, where the Frobenius argument used here is unavailable. The restriction on $q$ in the information-set upper bound also remains substantial for large $k$.

\end{document}